\documentclass[10pt,twoside]{article}
\usepackage{graphicx}
\usepackage{amsmath}
\usepackage{amssymb}

\usepackage{Latex-document}

\title{\bf On Some Inequalities for Gaussian \vskip -2mm Measures\vskip 6mm}

\author{R. Lata{\l}a\vspace*{-0.5cm}\thanks{Institute of Mathematics,
Warsaw University, Banacha 2, 02-097 Warszawa, Poland. E-mail:
rlatala@mimuw.edu.pl}}

\date{\vspace{-8mm}}

\def\er{\mathbb{R}}
\def\vol{\mathrm{vol}}

\begin{document}

\maketitle

\thispagestyle{first} \setcounter{page}{813}

\begin{abstract}

\vskip 3mm

We review several inequalities concerning Gaussian measures -
isoperimetric inequality, Ehrhard's inequality, Bobkov's inequality,
S-inequality and correlation conjecture.

\vskip 4.5mm

\noindent {\bf 2000 Mathematics Subject Classification:} 60E15,
60G15, 28C20, 26D15.

\noindent {\bf Keywords and Phrases:} Gaussian measure, Isoperimetry,
Ehrhard's inequality, Convex bodies, Correlation.
\end{abstract}

\vskip 12mm

\section{Introduction} \label{section 1}\setzero
\vskip-5mm \hspace{5mm }

Gaussian random variables and processes always played a central role in the
probability theory and statistics. The modern theory of Gaussian measures
combines methods from probability theory, analysis, geometry and topology
and is closely connected with diverse applications in functional analysis,
statistical physics, quantum field theory, financial mathematics and
other areas. Some examples of applications of Gaussian measures can be found in
monographs \cite{Bg,Le2, LT} and \cite{Lif}.

In this note we present several inequalities of geometric nature for
Gaussian measures. All of them have elementary formulations, but nevertheless
yield many important and nontrivial consequences. We begin in section 2 with
the already classical Gaussian isoperimetric inequality that inspired in the
70's and 80's the vigorous development of concentration inequalities
and their applications in the geometry and local theory of Banach spaces (cf.\
\cite{Le1,MS,Ta}).
In the sequel we review several more recent results and finish in section 6
with the discussion of the Gaussian correlation conjecture that remains
unsolved more than 30 years.

A probability measure $\mu$ on a real separable Banach space $F$ is
called {\it Gaussian} if for every functional $x^{*}\in F^{*}$ the
induced measure
$\mu \circ(x^{*})^{-1}$ is a one-dimensional Gaussian measure ${\cal
N}(a,\sigma^{2})$ for some $a=a(x^*)\in \er$ and $\sigma=\sigma(x^*)\geq 0$.
Throughout this note we only consider centered Gaussian measures that is the
measures such that $a(x^*)=0$ for all $x^*\in F^*$. A random vector with
values in $F$ is said to be Gaussian if its distribution is Gaussian.
Every centered Gaussian
measure on $\er^{n}$ is a linear image of the canonical Gaussian measure
$\gamma_{n}$, that is the measure on $\er^{n}$ with the density
$d\gamma_{n}(x)=(2\pi)^{-n/2} \exp(-|x|^{2}/2)dx$, where
$|x|=\sqrt{\sum_{i=1}^{n}x_{i}^{2}}$. Infinite dimensional Gaussian measures can
be effectively approximated by finite dimensional ones using the following
series representation (cf.\ \cite[Proposition 4.2]{Le2}): If
$\mu$ is a centered Gaussian measure on $F$ and $g_{1},g_{2},\ldots$ are
independent ${\cal N}(0,1)$ random variables then there exist vectors
$x_{1},x_{2},\ldots$ in $F$ such that the series
$X=\sum_{i=1}^{\infty}x_{i}g_{i}$ is convergent almost surely and in every
$L^{p}$, $0<p<\infty$, and is distributed as $\mu$.

We will denote by $\Phi$ the distribution function of the standard normal
${\cal N}(0,1)$ r.v., that is
\[\Phi(x)=\gamma_{1}(-\infty,x)=\frac{1}{\sqrt{2\pi}}\int_{-\infty}^{x}
e^{-y^{2}/2}dy, \ -\infty\leq x\leq \infty.\]
For two sets $A,B$ in a Banach space $F$ and $t\in \er$ we will write
$tA=\{tx:x\in A\}$ and $A+B=\{x+y:x\in A,y\in B\}$. A set $A$ in $F$ is said
to be {\it symmetric} if $-A=A$.

Many results presented in this note can be generalized to the more general
case of Radon Gaussian measures on locally convex spaces. For precise
definitions see \cite{Bg} or \cite{Bo3}.

\section{Gaussian isoperimetry} \label{section 2}
\setzero\vskip-5mm \hspace{5mm }

For a Borel set $A$ in $\er^{n}$ and $t>0$ let $A_{t}=A+tB_{2}^{n}=\{x\in
\er^{n}: |x-a|<t \mbox{ for some } a\in A\}$ be the open {\it t-enlargement} of
$A$, where $B_{2}^{n}$ denotes the open unit Euclidean ball in $\er^{n}$. The
classical isoperimetric inequality for the Lebesgue measure
states that if ${\rm vol}_{n}(A)={\rm vol}_{n}(rB_{2}^{n})$ then
${\rm vol}_{n}(A_{t})\geq {\rm vol}_{n}((r+t)B_{2}^{n})$ for $t>0$. In the early
70's C. Borell \cite{Bo2} and V.N. Sudakov and B.S. Tsirel'son \cite{ST} proved
independently the isoperimetric property of Gaussian measures.

{\bf Theorem 2.1} \it Let $A$ be a Borel set in $\er^{n}$ and let $H$ be
an affine halfspace such that $\gamma_{n}(A)=\gamma_{n}(H)=\Phi(a)$ for some
$a\in \er$. Then
\begin{equation}
\label{isop}
\gamma_{n}(A_{t})\geq \gamma_{n}(H_{t})=\Phi(a+t) \mbox{ for all } t\geq 0.
\end{equation}\rm

Theorem 2.1 has an equivalent differential analog. To state it let us define
for a measure $\mu$ on $\er^{n}$ and any Borel set $A$ the {\it boundary
$\mu$-measure} of $A$ by the formula
\[\mu^{+}(A)=\liminf_{t\rightarrow 0+}\frac{\mu(A_{t})-\mu(A)}{t}.\]
Moreover let $\varphi(x)=\Phi'(x)=(2\pi)^{-1/2}\exp(-x^{2}/2)$ and let
\[I(t)=\varphi\circ\Phi^{-1}(t), \ t\in [0,1]\]
be the {\it Gaussian isoperimetric function}.

The equivalent form of Theorem 2.1 is that for all Borel sets $A$ in
$\er^{n}$ \begin{equation}
\label{diffisop}
\gamma_{n}^{+}(A)\geq I(\gamma_{n}(A)).
\end{equation}
The equality in (\ref{diffisop})
holds for any affine halfspace.

For a probability measure $\mu$ on $\er^{n}$ we may define the {\it
isoperimetric function} of $\mu$ by
\[{\rm Is}(\mu)(p)=\inf\{\mu^{+}(A):\mu(A)=p\}, \ 0\leq p\leq 1.\]
Only few cases are known when one can determine exactly ${\rm Is}(\mu)$. For
Gaussian measures (\ref{diffisop}) states that ${\rm Is}(\gamma_{n})=I$.

Let us finish section 2 by an example of application of (\ref{isop}) (see
\cite[Lemma 3.1]{LT}).

{\bf Corollary 2.2} \it Let $X$ be a centered Gaussian random vector in a
separable Banach space $(F,\|\cdot\|)$. Then for any $t>0$
\[{\mathbf P}(|\|X\|-{\rm Med}(\|X\|)|\geq t)\leq 2(1-\Phi(\frac{t}{\sigma}))
\leq e^{-t^{2}/2\sigma^{2}},\]
where
\[\sigma=\sup\{\sqrt{{\mathbf E}(x^*(X))^{2}}:x^{*}\in F^{*},\|x^*\|\leq
1\}.\]\rm

\section{Ehrhard's inequality} \label{section 3}
\setzero\vskip-5mm \hspace{5mm}

It is well known that the classical isoperimetric inequality for the
Lebesgue measure in $\er^{n}$ follows by the Brunn-Minkowski
inequality (cf.\ \cite{Os}), which states that for any Borel sets $A$ and $B$ in
$\er^{n}$ \[\vol_{n}(\lambda  A+ (1-\lambda)B))\geq (\vol_{n}(A))^{\lambda}
  (\vol_{n}(B))^{1-\lambda} \mbox{ for }\lambda\in [0,1].\]
Gaussian measures
satisfy the similar log-concavity property, that is the inequality
\begin{equation}
\label{logconc}
\ln(\mu(\lambda A+ (1-\lambda)B))\geq \lambda
  \ln(\mu(A))+(1-\lambda)\ln(\mu(B)),\ \lambda\in[0,1]
\end{equation}
holds for any Gaussian measure $\mu$ on a separable Banach space $F$ and any
Borel sets $A$ and $B$ in $F$ (cf.\ \cite{Bo1}). However the log-concavity of
the measure does not imply the Gaussian isoperimetry.

In the early 80's A.\ Ehrhard \cite{Eh} gave a different proof of the
isoperimetric inequality (\ref{isop}) using a Gaussian symmetrization procedure
similar to the Steiner symmetrization. With the same symmetrization tool Ehrhard
established a new  Brunn-Minkowski type inequality, stronger than
(\ref{logconc}), however only for convex sets.

{\bf Theorem 3.1}(Ehrhard's inequality) \it If $\mu$ is a centered Gaussian
measure on a separable Banach space $F$ and $A$, $B$ are Borel sets in $F$, with
at least one of them convex, then
\begin{equation}
\label{ehr}
\Phi^{-1}(\mu(\lambda A+ (1-\lambda)B))\geq \lambda
  \Phi^{-1}(\mu(A))+(1-\lambda)\Phi^{-1}(\mu(B)) \mbox{ for }\lambda\in [0,1].
\end{equation} \rm

For both sets $A$ and $B$ convex Ehrhard's inequality was proved in
\cite{Eh}. The generalization to the case when only one of the sets is convex
was established in \cite{La}.

It is not hard to see that Theorem 3.1 implies the isoperimetric inequality
(\ref{isop}). Indeed we have for any Borel set $A$ in $\er^{n}$
\begin{eqnarray*}
\lefteqn{\Phi^{-1}(\gamma_{n}(A_{t}))=\Phi^{-1}(\gamma_{n}(\lambda
(\lambda^{-1}A)+ (1-\lambda)((1-\lambda)^{-1}tB_{2}^{n})))}\\
& &\geq
\lambda\Phi^{-1}(\gamma_{n}(\lambda^{-1}A))+
(1-\lambda)\Phi^{-1}(\gamma_{n}((1-\lambda)^{-1}tB_{2}^{n}))
\overset{\lambda\rightarrow 1-}{\longrightarrow}
\Phi^{-1}(\gamma_{n}(A))+t.
\end{eqnarray*}

{\bf Conjecture 3.1} \it Inequality (\ref{ehr}) holds for any Borel sets
in $F$. \rm

Ehrhard's symmetrization procedure enables us to reduce Conjecture 3.1 to the
case $F=\er$ and $\mu=\gamma_{1}$. We may also assume that $A$ and $B$ are
finite unions of intervals. At the moment the conjecture is known to hold
when $A$ is a union of at most 3 intervals.

Ehrhard's inequality has the following Prekopa-Leindler type
functional version. Suppose that $\lambda\in (0,1)$ and
$f,g,h:\er^{n}\rightarrow [0,1]$ are such that
\[\forall_{x,y\in \er^{n}}\ \Phi^{-1}(h(\lambda x+(1-\lambda)y))\geq
 \lambda\Phi^{-1}(f(x))+(1-\lambda)\Phi^{-1}(g(y))\]
then
\begin{equation}
\label{funehr}
\Phi^{-1}(\int_{\er^{n}}hd\gamma_{n})\geq \lambda
\Phi^{-1}(\int_{\er^{n}}fd\gamma_{n})+(1-\lambda)
\Phi^{-1}(\int_{\er^{n}}gd\gamma_{n}).
\end{equation}
We use here the convention $\Phi^{-1}(0)=-\infty,\Phi^{-1}(1)=\infty$
and $-\infty+\infty=-\infty$. At the moment the above functional inequality
is known to hold under the additional assumption that at least one of the
functions $\Phi^{-1}(f), \Phi^{-1}(g)$ is convex. When one
takes $f=1_{A}$, $g=1_{B}$ and $h=1_{\lambda A+(1-\lambda)B}$ the inequality
(\ref{funehr}) immediately implies (\ref{ehr}). On the other hand if we put
$A=\{(x,y)\in \er^{n}\times\er: y\leq \Phi^{-1}(f(x))\}$ and $B=\{(x,y)\in
\er^{n}\times\er: y\leq \Phi^{-1}(g(x))\}$ then $\lambda
A+(1-\lambda)B \subset \{(x,y)\in \er^{n}\times\er: y\leq \Phi^{-1}(h(x))\}$,
so Ehrhard's inequality in $\er^{n+1}$ implies (\ref{funehr})
in $\er^{n}$. It is easy to show the inductive step in the proof of
$(\ref{funehr})$. Unfortunately the case $n=1$ in the functional inequality
seems to be much more complicated than the case $\mu=\gamma_{1}$ in Ehrhard's
inequality.

\section{Bobkov's inequality} \label{section 4}
\setzero\vskip-5mm \hspace{5mm}

Isoperimetric inequality for the Lebesgue measure has an equivalent analytic
form - the Sobolev inequality (cf.\ \cite{Os}). L. Gross
\cite{Gr}
showed that the Gaussian measures $\gamma_{n}$ satisfy the logarithmic
Sobolev inequality
\begin{equation}
\label{logsob}
 \int_{\er^{n}}g^{2}\log g^{2}d\gamma_{n}-\int_{\er^{n}}g^{2}d\gamma_{n}
 \log(\int_{\er^{n}}g^{2}d\gamma_{n})\leq
 2 \int_{\er^{n}}|\nabla g|^{2}d\gamma_{n}
\end{equation}
for all smooth functions $g:\er^{n}\rightarrow \er$. Using the so-called
Herbst argument one can show (cf.\ \cite[Sect.\ 5.1]{Le1}) that
(\ref{logsob}) implies the concentration inequality
\[\gamma_{n}(\{h\geq \int_{\er^{n}}hd\gamma_{n}+t\})\leq e^{-t^{2}/2},\
  t\geq 0\]
valid for all Lipschitz functions $h:\er^{n}\rightarrow \er$ with the Lipschitz
seminorm $\|h\|_{\rm Lip}=\sup\{|h(x)-h(y)|:x,y\in \er^{n}\}\leq 1$.
However the logarithmic Sobolev inequality does not imply the isoperimetric
inequality.

The formulation of the functional form of Gaussian isoperimetry
was given by S.G.\ Bobkov \cite{Bb}.

{\bf Theorem 4.1} \it For any locally Lipschitz function $f:\er^{n}\rightarrow
[0,1]$ and $\mu=\gamma_{n}$ we have
\begin{equation}
\label{bobin}
I(\int_{\er^{n}}f d\mu)\leq \int_{\er^{n}}
 \sqrt{I(f)^{2}+|\nabla f|^{2}}d\mu.
\end{equation}\rm

Theorem 4.1 easily implies the isoperimetric inequality (\ref{diffisop}) by
approximating the indicator function $I_{A}$ by Lipschitz functions. On the
other hand if we apply (\ref{diffisop}) to the set $A=\{(x,y)\in
\er^{n}\times\er: \Phi(y)<f(x)\}$ in $\er^{n+1}$ we get (\ref{bobin}). It is
also not hard to derive the logarithmic Sobolev inequality (\ref{logsob}) as a
limit case of Bobkov's inequality (cf. \cite{BM}): one should use (\ref{bobin})
for $f=\varepsilon g^{2}$ (with $g$ bounded) and let $\varepsilon$ tend to 0
($I(t)\sim t\sqrt{2\log (1/t)}$ as $t\rightarrow 0+$).

The crucial point of the inequality (\ref{bobin}) is its tensorization
property. To state it precisely let us say that a measure $\mu$ on $\er^{n}$
{\it satisfies Bobkov's inequality} if the inequality (\ref{bobin}) holds
for all locally Lipschitz functions $f:\er^{n}\rightarrow [0,1]$. Easy
argument shows that if $\mu_{i}$ are measures on $\er^{n_{i}}$, $i=1,2$,
that satisfy Bobkov's inequality then the measure $\mu_{1}\otimes\mu_{2}$
also satisfies Bobkov's inequality.

The inequality (\ref{bobin}) was proved by Bobkov in an elementary way,
based on the following "two-point" inequality:
\begin{equation}
\label{twop}
I(\frac{a+b}{2})\leq
\frac{1}{2}\sqrt{I(a)^{2}+(\frac{a-b}{2})^{2}}+
\frac{1}{2}\sqrt{I(b)^{2}+(\frac{a-b}{2})^{2}}
\end{equation}
valid for all $a,b\in [0,1]$.
In fact the inequality (\ref{twop}) is equivalent to Bobkov's inequality for
$\mu=\frac{1}{2}\delta_{-1}+ \frac{1}{2}\delta_{1}$ and the discrete gradient
instead of $\nabla f$. Using the tensorization property and the central limit
theorem Bobkov deduces (in the similar way as Gross in his proof of
(\ref{logsob})) (\ref{bobin}) from (\ref{twop}).

Using the co-area formula and Theorem 4.1 F.\ Barthe and M.\ Maurey
\cite{BM} gave
interesting characterization of all absolutely continuous measures
that satisfy Bobkov's inequality.

{\bf Theorem 4.2} \it Let $c>0$ and $\mu$ be a Borel probability measure on
the Riemannian manifold $M$, absolutely continuous with respect to the
Riemannian volume. Then the following properties are equivalent\\
(i) For every measurable $A\subset M$, $\mu^{+}(A)\geq cI(\mu(A))$;\\
(ii) For every locally Lipschitz function $f:M\rightarrow [0,1]$
\[I(\int_{M}f d\mu)\leq \int_{M}
 \sqrt{I(f)^{2}+\frac{1}{c^{2}}|\nabla f|^{2}}d\mu.\]\rm

Theorem 4.2 together with the tensorization property shows that if
${\rm Is}(\mu_{i})\geq cI$, $i=1,2\ldots$, then also
${\rm Is}(\mu_{1}\otimes\ldots\otimes\mu_{n})\geq cI$. In general it is
not known how to estimate ${\rm Is}(\mu_{1}\otimes\ldots\otimes\mu_{n})$
in terms of ${\rm Is}(\mu_{i})$ even in the case when all $\mu_{i}$'s are
equal (another important special case of this problem
was solved in \cite{BH}) .

\section{S-inequality} \label{section 5}
\setzero\vskip-5mm \hspace{5mm}

In many problems arising in probability in Banach spaces one needs to estimate
the measure of balls in some Banach space $F$. In particular one may ask
what is the slowest possible grow of the Gaussian measure of balls in $F$ or
more general of some fixed convex symmetric closed set under dilations. The next
theorem, proved by R. Lata{\l}a and K. Oleszkiewicz \cite{LO}, gives the
positive answer to the conjecture posed in an unpublished manuscript of L.\ A.\
Shepp (1969).

{\bf Theorem 5.1}(S-inequality) \it Let $\mu$ be a centered Gaussian measure
on a separable Banach space $F$. If $A$ is a symmetric, convex, closed subset
of $F$ and $P\subset F$ is a symmetric strip, that is
$P=\{x\in F:|x^{*}|\leq 1\}$ for some $x^*\in F^*$, such that $\mu(A)=\mu(P)$
then \[\mu(tA)\geq \mu(tP) \mbox{ for } t\geq 1\]
and
\[\mu(tA)\leq \mu(tP) \mbox{ for } 0\leq t\leq 1.\] \rm

A simple approximation argument shows that it is enough to prove Theorem 5.1
for $F=\er^{n}$ and $\mu=\gamma_{n}$. The case $n\leq 3$ was solved by
V.N.\ Sudakov
and V.A.\ Zalgaller \cite{SZ}. Under the additional assumptions of symmetry of
$A$ in $\er^{n}$ with respect to each coordinate, Theorem 5.1 was proved
by S.\ Kwapie\'n and J.\ Sawa \cite{KS}.

S-inequality can be equivalently expressed as
\[\Psi^{-1}(\mu(tA))\geq t\Psi^{-1}(\mu(A)) \mbox{ for } t\geq 1,\]
where $\Psi^{-1}$ denotes the inverse of
\[\Psi(x)=\gamma_{1}(-x,x)=\frac{1}{\sqrt{2\pi}}\int_{-x}^{x}e^{-y^{2}/2}dy.\]

The crucial tool in the proof of S-inequality is the new modified isoperimetric
inequality. Let us first define for a convex symmetric set $A$ in $\er^{n}$
\[w(A)=2\sup\{r: B(0,r)\subset A\}.\]
It is easy to see that for a symmetric strip $P$, $w(P)$ is equal to the width
of $P$ and for a symmetric convex set $A$
\begin{equation}
\label{defwidth}
w(A)=\inf\{w(P): A\subset P, P \mbox{ is a symmetric strip in }\er^{n}\}.
\end{equation}
Thus $w(A)$ can be considered as the width of the set $A$. The following
isoperimetric-type theorem holds true.

{\bf Theorem 5.2} \it If $\gamma_{n}(A)=\gamma_{n}(P)$, where P is a symmetric
strip and $A$ is a convex symmetric set in $\er^{n}$, then
\begin{equation}
\label{modisop}
w(A)\gamma_{n}^{+}(A)\geq w(P)\gamma_{n}^{+}(P).
\end{equation}\rm

The main advantage of the inequality (\ref{modisop})  is that one may
apply here the symmetrization procedure and reduce Theorem 5.2 to the similar
statement for 2-dimensional convex sets symmetric with respect to some axis.

It is not hard to see that Theorem 5.2 implies Theorem 5.1. Indeed, let us
define for any measurable set $B$ in $\er^{n}$, $\gamma_{B}(t)=\gamma_{n}(tB)$
for $t>0$. Taking the derivatives of both sides of the inequalities in Theorem
5.1 one can see that it is enough to show
\begin{equation}
\label{diffsin}
\gamma_{n}(A)=\gamma_{n}(P) \Rightarrow \gamma_{A}'(1)\geq \gamma_{P}'(1)
\end{equation}
for any symmetric convex closed set $A$ and a symmetric strip
$P=\{|x_{1}|\leq p\}$. Let $w=w(A)$, so $B(0,w)\subset A$. Then for $t>1$ and
$x\in A$ we have $B(t^{-1}x,(t-1)w/t)=t^{-1}x+(1-t^{-1})B(0,w) \subset A$, so
$B(x, (t-1)w)\subset tA$. Hence $A_{(t-1)w}\subset tA$ and
\[\gamma_{A}'(1)\geq w\gamma^{+}_{n}(A)=w(A)\gamma^{+}_{n}(A).\]
However for the strip $P$
\[\gamma_{P}'(1)=\sqrt{\frac{2}{\pi}}pe^{-p^{2}/2}=
 w(P)\gamma^{+}_{n}(P)\]
and the inequality (\ref{diffsin}) follows by Theorem 5.2.

It is not clear if the convexity assumption for the set $A$ in Theorem 5.2
is necessary (obviously $w(A)$ for nonconvex symmetric sets $A$ should be
defined by (\ref{defwidth})). One may also ask if the symmetry assumption
can be released (with the suitable modification of the definition of the
width for nonsymmetric sets). Also functional versions of Theorems 5.1
and 5.2 are not known.

As was noticed by S. Szarek S-inequality implies the best constants in
comparison of moments of Gaussian vectors (cf.\ \cite{LO}).

{\bf Corollary 5.3} \it If $X$ is a centered Gaussian vector in a separable
Banach space $(F,\|\cdot\|)$ then
\[({\mathbf E}\|X\|^{p})^{1/p}\leq \frac{c_{p}}{c_{q}}
  ({\mathbf E}\|X\|^{q})^{1/q}
 \mbox{ for any } p\geq q\geq 0,\]
where
\[c_{p}=({\mathbf E}|g_{1}|^{p})^{1/p}=\sqrt{2}(\frac{1}{\sqrt{\pi}}
  \Gamma(\frac{p+1}{2}))^{1/p}.\] \rm

Another interesting problem connected with the S-inequality was recently posed
by W.\ Banaszczyk (private communication): Is it true that under the
assumptions of Theorem 5.1
\begin{equation}
\label{bana}
\mu(s^{\lambda}t^{1-\lambda}A)\geq \mu(sA)^{\lambda}\mu(tA)^{1-\lambda},\
\lambda\in[0,1]
\end{equation}
for any closed convex symmetric set $A$ in $F$ and $s,t>0$? Combining the facts
that the function $\Phi^{-1}(\mu(tA))$ is concave (Theorem 3.1) and the function
$\frac{1}{t}\Psi^{-1}(\mu(tA))$ is nondecreasing (Theorem 5.1) one can show
that (\ref{bana}) holds if $\mu(sA),\mu(tA)\geq c$, where $c<0.85$ is some
absolute constant.

It is of interest if Theorem 5.1 can be extended to the more general class
of measures. The following conjecture seems reasonable.

{\bf Conjecture 5.1} \it Let $\nu$ be a rotationally invariant measure
on $\er^{n}$, absolutely continuous with respect to the Lebesgue measure
with the density of the form $f(|x|)$ for
some nondecreasing function $f:\er_{+}\rightarrow [0,\infty)$. Then for
any convex symmetric set $A$ in $\er^{n}$ and any symmetric strip $P$ in
$\er^{n}$ such that $\nu(A)=\nu(P)$ the inequality $\nu(\lambda A)\geq
\nu(\lambda P)$ is satisfied for $\lambda\geq 1$.\rm

To show Conjecture 5.1 it is enough to establish the following conjecture
concerning the volumes of the convex hulls of symmetric sets on the
$n-1$-dimensional unit sphere $S^{n-1}$.

{\bf Conjecture 5.2} \it Let $\sigma_{n-1}$ be a Haar measure on $S^{n-1}$,
$A$ be a symmetric subset of $S^{n-1}$ and $P=\{x\in S^{n-1}: |x_{1}|\leq
t\}$ be a symmetric strip on $S^{n-1}$ such that $\sigma_{n-1}(A)=
\sigma_{n-1}(P)$, then $\vol_{n}(\mathrm{conv}( A))\geq
\vol_{n}(\mathrm{conv}(P))$.\rm

It is known that both conjectures hold for $n\leq 3$ (cf.\ \cite{SZ}).

\section{Correlation conjecture} \label{section 6}
\setzero\vskip-5mm \hspace{5mm}

The following conjecture is an object of intensive efforts of many
probabilists since more then 30 years.

{\bf Conjecture 6.1} \it If $\mu$ is a centered Gaussian measure on a
separable Banach space $F$ then
\begin{equation}
\label{corr}
\mu(A\cap B)\geq \mu(A)\mu(B)
\end{equation}
for all convex symmetric sets $A,B$ in $F$.\rm

Various equivalent formulations of
Conjecture 6.1 and  history of the problem can be found in \cite{SSZ}.
Standard approximation argument shows that it is enough to show
(\ref{corr}) for $F=\er^{n}$ and $\mu=\gamma_{n}$. For $n=2$ the solution was
given by L. Pitt \cite{Pi}, for $n\geq 3$ the conjecture remains unsettled,
but a variety of special results are known. Borell \cite{Bo4} established
(\ref{corr}) for sets $A,B$  in a certain class of (not necessary
convex) sets in $\er^{n}$, which for $n=2$ includes all symmetric sets. A
special case of (\ref{corr}), when one of the sets $A,B$ is a symmetric strip of
the form $\{x\in F: |x^{*}(x)|\leq 1\}$ for some $x^*\in F^{*}$, was proved
independently by C.\ G.\ Khatri \cite{Kh} and Z.\ \v{S}id\'ak \cite{Si}
(see \cite{GE} for an extension to elliptically
contoured distributions and \cite{SW} for the case when one of the
sets is a nonsymmetric strip). Recently,
the Khatri-\v{S}id\'ak result has been generalized by G. Harg\'e \cite{Ha} to
the case when one of the sets is a symmetric ellipsoid.

{\bf Theorem 6.1} \it If $\mu$ is a centered Gaussian measure on $\er^{n}$, $A$
is a symmetric convex set in $\er^{n}$ and $B$ is a symmetric ellipsoid, that
is the set of the form $B=\{x\in \er^{n}: \langle Cx,x\rangle\leq 1\}$ for
some symmetric nonnegative matrix $C$, then
\[\mu(A\cap B)\geq \mu(A)\mu(B).\]\rm

The following weaker form of (\ref{corr})
\[\mu(A\cap B)\geq \mu(\lambda A)\mu(\sqrt{1-\lambda^{2}}B), \
0\leq \lambda\leq 1\]
was established for $\lambda=\frac{1}{\sqrt{2}}$ in \cite{SSZ} and for
general $\lambda$ in \cite{Li}. The Khatri-\v{S}id\'ak result and the above
inequality turn out to be very useful in the study of the so-called small ball
probabilities for Gaussian processes (see \cite{LiS} for a survey of results in
this direction).

The correlation conjecture has the following functional form:
\begin{equation}
\label{funcorr}
\int fg d\mu\geq \int f d\mu\int g d\mu
\end{equation}
for all nonnegative even functions $f,g$ such that the sets $\{f\geq t\}$ and
$\{g\geq t\}$ are convex for all $t\geq 0$. Y.\ Hu \cite{Hu} showed that
the inequality (\ref{funcorr}) (that we would like to have for log-concave
functions) is valid for even convex functions $f,g\in L^{2}(F,\mu)$.

\label{lastpage}


\begin{thebibliography}{aa}
\bibitem{BM} F. Barthe, B. Maurey, Some remarks on isoperimetry of Gaussian
type, {\it Ann. Inst. H. Poincar\'e Probab. Statist.} 36 (2000), 419--434.
\bibitem{Bb} S.G. Bobkov, An isoperimetric inequality on the discrete cube, and
an elementary proof of the isoperimetric inequality in Gauss space, {\it
Ann. Probab.} 25 (1997), 206--214.
\bibitem{BH} S.G. Bobkov, C. Houdr\'e, Isoperimetric constants for product
probability measures, {\it Ann. Probab.} 25 (1997), 184--205.
\bibitem{Bg} V.I. Bogachev, {\it Gaussian Measures}, American Mathematical
Society, Providence, RI, 1998.
\bibitem{Bo1} C. Borell, Convex measures on locally convex spaces, {\it
Ark. Mat.} 12 (1974), 239--252.
\bibitem{Bo2} C. Borell, The Brunn-Minkowski inequality in Gauss space,
{\it Invent. Math.}, 30 (1975), 207--216.
\bibitem{Bo3} C. Borell, Gaussian Radon measures on locally convex spaces,
{\it Math. Scand.} 38 (1976), 265--284.
\bibitem{Bo4} C. Borell, A Gaussian correlation inequality for certain bodies
in $R^{n}$, {\it Math. Ann.} 256 (1981), 569--573.
\bibitem{Eh} A. Ehrhard, Sym\'etrisation dans l'espace de Gauss,
{\it Math. Scand.}, 53 (1983), 281--301.
\bibitem{Gr} L. Gross,  Logaritmic Sobolev inequalities, {\it Amer. J. Math.}
97 (1975), 1061--1083.
\bibitem{GE} S. Das Gupta, M.L. Eaton, I. Olkin, M. Perlman, L.J. Savage,
M. Sobel, Inequalities on the probability content of convex regions
for elliptically contoured distributions, {\it Proc. Sixth Berkeley
Symp. Math. Statist. Prob.} vol. II, 241--264, Univ. California Press,
Berkeley, 1972.
\bibitem{Ha} G. Harg\'e, A particular case of correlation inequality for
the Gaussian measure, {\it Ann. Probab.} 27 (1999), 1939--1951.
\bibitem{Hu} Y. Hu, It\^{o}-Wiener chaos expansion with exact residual and
correlation, variance inequalities, {\it J. Theoret. Probab.} 10 (1997),
835--848.
\bibitem{Kh} C.G. Khatri, On certain inequalities for normal distributions
and their applications to simultaneous confidence bounds, {\it Ann. Math.
Stat.} 38 (1967), 1853--1867.
\bibitem{KS} S. Kwapie\'n, J. Sawa, On some conjecture concerning
Gaussian measures of dilatations of convex symmetric sets, {\it Studia Math.}
105 (1993), 173--187.
\bibitem{La} R. Lata{\l}a, A note on the Ehrhard inequality, {\it Studia
Math.} 118 (1996), 169--174.
\bibitem{LO} R. Lata{\l}a, K. Oleszkiewicz, Gaussian measures of dilatations
of convex symmetric sets, {\it Ann. Probab.} 27 (1999), 1922--1938.
\bibitem{Le2} M. Ledoux, Isoperimetry and Gaussian Analysis, {\it Lectures
on probability theory and statistics (Saint-Flour, 1994)}, 165--294,
Lecture Notes in Math. 1648, Springer, Berlin, 1996.
\bibitem{Le1} M. Ledoux, {\it The concentration of measure phenomenon},
American Mathematical Society, Providence, RI, 2001.
\bibitem{LT} M. Ledoux, M. Talagrand, {\it Probability on Banach Spaces.
Isoperimetry and processes}, Springer-Verlag, Berlin, 1991.
\bibitem{Li} W.V. Li, A Gaussian correlation inequality and its applications
to small ball probabilities, {\it Electron. Comm. Probab.} 4 (1999), 111-118.
\bibitem{LiS} W.V. Li, Q.M. Shao, Gaussian processes: inequalities, small
ball probabilities and applications, {\it Stochastic Processes: Theory
and Methods}, Handbook of Statistics vol. 19, 533--597, Elsevier,
Amsterdam 2001.
\bibitem{Lif} M.A. Lifshits, {\it Gaussian random functions}, Kluwer
Academic Publications, Dordrecht, 1995.
\bibitem{MS} V.D. Milman, G. Schechtman, {\it Asymptotic theory of
finite-dimensional normed spaces}, Lecture Notes in Math. 1200, Springer-Verlag,
Berlin, 1986.
\bibitem{Os} R. Osserman, The isoperimetric inequality, {\it Bull.
Amer. Math. Soc.} 84 (1978), 1182--1238.
\bibitem{Pi} L. Pitt, A Gaussian correlation inequality for symmetric convex
sets, {\it Ann. Probability}, 5 (1977), 470--474.
\bibitem{SSZ} G. Schechtman, T. Schlumprecht, J. Zinn, On the Gaussian
measure of intersection, {\it Ann. Probab.} 26 (1998), 346--357.
\bibitem{Si} Z. \v{S}id\'ak, Rectangular confidence regions for the means of
multivariate normal distributions, {\it J. Amer. Statist. Assoc.} 62 (1967),
626--633.
\bibitem{ST} V.N. Sudakov, B.S. Tsirel'son, Extremal propertiesof half-spaces
for spherically invariant measures (in Russian), {\it Zap. Nauchn. Sem.
L.O.M.I.} 41 (1974), 14--24.
\bibitem{SZ} V.N. Sudakov, V.A. Zalgaller, Some problems on centrally
symmetric convex bodies (in Russian) {\it Zap. Nauchn. Sem. L.O.M.I.} 45
(1974), 75--82.
\bibitem{SW} S. Szarek, E. Werner, A nonsymmetric correlation inequality
for Gaussian measure, {\it J. Multivariate Anal.} 68 (1999), 193--211.
\bibitem{Ta} M. Talagrand, Concentration of measure and isoperimetric
inequalities in product spaces, {\it IHES Publ. Math.} 81 (1995), 73--205.
\end{thebibliography}
\end{document}